\documentclass[letterpaper,10pt,oneside,onecolumn,reqno]{amsart}

\usepackage{amsmath, amssymb}
\usepackage{url}
\usepackage[all]{xy}
\usepackage[left=.95in,top=.95in,right=.95in,bottom=.95in]{geometry}

\newcommand{\B}{\mathcal B}

\newcommand{\CC}{\mathcal C}
\newcommand{\EE}{\mathbb E}

\newcommand{\PP}{\mathbb P}

\newcommand{\R}{\mathbb R}

\newcommand{\supp}{\operatorname{supp}}
\renewcommand{\span}{\operatorname{span}}
\newcommand{\D}{\mathrm{d}}				
\newcommand{\sD}{\, \mathrm{d}}		
\newcommand{\oo}{\infty}

\numberwithin{equation}{section}
\theoremstyle{definition}
\newtheorem{env_thmmain}{Theorem}
\newtheorem{env_cormain}[env_thmmain]{Corollary}
\newtheorem{env_thm}{Theorem}[section]
\newtheorem{env_lem}[env_thm]{Lemma}

\newtheorem{env_def}[env_thm]{Definition}

\title{Continuous Disintegrations of Gaussian Processes}	
\author{T. LaGatta}
\date{Fall 2010}

\begin{document}

	\begin{abstract}
		The goal of this paper is to understand the conditional law of a stochastic process once it has been observed over an interval.  To make this precise, we introduce the notion of a continuous disintegration:  a regular conditional probability measure which varies continuously in the conditioned parameter.  The conditioning is infinite-dimensional in character, which leads us to consider the general case of probability measures in Banach spaces.  Our main result is that for a certain quantity $M$ based on the covariance structure, $M < \oo$ is a necessary and sufficient condition for a Gaussian measure to have a continuous disintegration.  The condition $M < \oo$ is quite reasonable: for the familiar case of stationary processes, $M = 1$.
	\end{abstract}

	\maketitle
	
	\section{Introduction}
	
	Consider a continuous Gaussian process $\xi_t$ on an interval $[0,T]$.  Let $S \le T$, and let $y(s)$ be a continuous function on the sub-interval $[0,S]$.  Suppose that we observe $\xi_s = y(s)$ for all $s \le S$.  This paper is a result of asking the following questions:  
		\begin{itemize}
			\item Is the conditional law $\PP^y := \PP\big( \cdot \big|\, \xi|_{[0,S]} = y \big)$ still Gaussian?
			\item Is there a sufficient condition so that the measures $\PP^y$ vary continuously in the parameter $y$?
		\end{itemize}
	The answer to both these questions is ``yes'', as illustrated by the following theorem.
	

	\begin{env_thmmain} \label{thm_stoch} 
		Let $\xi_t$ be a Gaussian process on $[0,T]$ with mean zero and covariance function $c$:
			$$\EE \xi_t = 0 \qquad \mathrm{and} \qquad \EE \xi_t \xi_s = c(t,s).$$
		Suppose that $\xi_t$ is almost-surely continuous.  Let $S \le T$, and suppose that
			\begin{equation} \label{stoch_M}
				M = \sup_{s \le S} \frac{\sup_{t \le T} |c(s,t)|}{\sup_{s' \le S} |c(s,s')|} < \oo. \end{equation}
		There exists a closed family $Y_0$ of functions such that with probability one, $\xi|_{[0,S]} \in Y_0$, and the regular conditional probability $\PP^y := \PP\big( \cdot \big|\, \xi|_{[0,S]} = y \big)$ is a well-defined Gaussian measure which varies (weakly) continuously in $y \in Y_0$. 

		If $\xi_t$ is a stationary process, then $M = 1$.
	\end{env_thmmain}
	
	The case that $\xi_t$ has mean function $\mu(t)$ is handled by applying this theorem to the mean-zero stochastic process $\xi_t - \mu(t)$.  This theorem is a special case of Theorem \ref{thm_RF}, which applies to continuous random fields (stochastic processes) defined over compact parameter spaces.  
	
	Since the conditioning is of a function over an entire interval, it is infinite-dimensional in character.  This leads us to consider the Banach spaces $X = C([0,T])$ and $Y = C([0,S])$ of continuous functions equipped with the $\sup$ norms, as well as the restriction map $\eta : X \to Y$.  The main result of the paper, Theorem \ref{thm_mainresult}, is simply the general form of the above theorem in the context of arbitrary Banach spaces.
	
	To explain what it means for a regular conditional probability measure to vary (weakly) continuously, we introduce the notion of a notion of a continuous disintegration of a probability measure $\PP$.
		
	\begin{env_def}		
	Let $X$ and $Y$ be complete metric spaces, with Borel $\sigma$-algebras $\B(X)$ and $\B(Y)$, and let $\PP$ be a Radon probability measure on $X$.  Let $\eta : X \to Y$ be a measurable function, and denote the push-forward measure of $\PP$ on $Y$ by $\PP_Y = \PP \circ \eta^{-1}$.  A \emph{disintegration} (or \emph{regular conditional probability}) of $\PP$ with respect to $\eta$ is a map $Y \times \B(X) \to \R$ (denoted by $(y,B) \mapsto \PP^y(B)$) such that:
		\begin{itemize}
			\item For all $y \in Y$, $\PP^y$ is a probability measure on $\B(X)$.
			\item For all $B \in \B(X)$, $y \mapsto \PP^y(B)$ is a measurable function of $y \in Y$.
			\item The measure $\PP^y$ is supported on the fiber of $y$.  i.e., for $\PP_Y$-almost every $y \in Y$, $\PP^y(\eta^{-1}(y)) = 1$, and
			\item For all integrable functions $f : X \to \R$, the disintegration equation holds:
				\begin{equation} \label{disint}
					\int_X f(x) \sD \PP(x) = \int_Y \int_X f(x) \sD \PP^y(x) \D \PP_Y(y). \end{equation}
		\end{itemize} 
	Suppose furthermore that $\eta$ is continuous and $Y_0$ is a closed subset of $Y$ of full $\PP_Y$-measure.  We say that $\PP^y$ is a \emph{continuous disintegration} given $Y_0$ provided
		\begin{equation} \label{contdis}
			\mbox{if $y_n \in Y_0$ and $y_n \to y$, then $\mathbb P^{y_n}$ converges weakly to $\mathbb P^y$.} \end{equation}
	\end{env_def}
	
	We remark that disintegration is typically a more general concept than regular conditional probability.  In the present work, we ignore the distinction and treat the terms as synonyms.
	

	The notion of a continuous disintegration is a new contribution to the literature, but disintegrations and regular conditional probabilities have been studied in wide generality.  For a nice overview of the topic, see the survey \cite{chang1997conditioning} by Chang and Pollard, or Sections 10.4 and 10.6 of Bogachev \cite{bogachev2007measure_vol2}.  The typical existence theorem can be found in \cite[Section 4.1.c]{durrett1996probability} or \cite[Theorem 1.4]{billingsley1968convergence}.  The recent paper  \cite{leao2004regular} contains some very general existence results.
	
	Our main result, Theorem \ref{thm_mainresult}, gives a sufficient condition for continuous disintegrations to exist for Gaussian probability measures on Banach spaces.  Theorem \ref{thm_Minfinity} demonstrates that this condition is also necessary.  Theorem \ref{thm_RF} is the application of the existence result to the important context of random fields.\newline
		
	
	Suppose that $X$ and $Y$ are finite-dimensional vector spaces, the map $\eta : X \to Y$ is linear, and the measure $\PP$ is Gaussian.  It is a simple exercise in linear algebra that the regular conditional probability measure $\PP^y = \PP(\cdot|\, \eta^{-1}(y))$ is Gaussian, and that the conditioned covariance matrix does not depend on the actual value $y \in Y$.  The conditional mean vector is easily seen to vary continuously in $y$.  Since $\PP^y$ is Gaussian and depends entirely on its mean and covariance matrix, it follows easily that $\PP^y$ is a continuous disintegration.
	
	Now suppose that the spaces $X$ and $Y$ are separable Banach spaces, the map $\eta : X \to Y$ is linear and continuous, and the measure $\PP$ is Gaussian.  Tarieladze and Vakhania \cite{tarieladze2007disintegration} show that $\PP$ admits a disintegration $\PP^y$ which is a Gaussian measure for all $y$.  Furthermore, when the push-forward measure $\PP_Y$ has finite-dimensional support in $Y$, it quickly follows from their Theorem 3.11a that $\PP^y$ is a continuous disintegration given $\supp \PP_Y$.  
	
	This fact is quite useful in applications, such as kriging in geosciences and hydrology \cite{krigingnotes}.  In this example, one models a quantity of interest, such as elevation, by a Gaussian random field defined on a domain in $\R^2$.  By conditioning the field at finitely many points based on empirical data, the field serves as a reasonable interpolation between the sampled points, with the randomness representing uncertainty.  The result of Tarieladze and Vakhania demonstrates that a conditioned Gaussian field is still Gaussian, and that its law varies continuously in the sampled values.\newline
	
	In this paper, we focus on the situation where $X$ and $Y$ are arbitrary Banach spaces, the map $\eta : X \to Y$ is linear and continuous, and $\PP$ is a Gaussian measure with mean zero and covariance operator $K$.  We need not worry about separability of $X$, as the structure theorem (Theorem \ref{thm_structure}) asserts that the Radon measure $\PP$ is supported on the separable subspace $\overline{KX^*}$.  The push-forward measure $\PP_Y$ has covariance operator $\eta K \eta^*$, and is supported on $Y_0 := \overline{\eta K \eta^* Y^*}$.  
	
	In Lemma \ref{inj}, we show that the map $\eta$ is injective when restricted to $K\eta^*Y^* \subseteq X$.  Consequently, the inverse map $\eta^{-1} : \eta K \eta^* Y^* \to X$ is well-defined.  Let $M$ denote the operator norm of $\eta^{-1}$, and suppose that $M < \oo$.  Then $\eta^{-1}$ extends to a continuous linear map $m : Y_0 \to X$.  In Lemma \ref{khat}, we use a Hilbert-space formalism to show that the operator $\hat K := K - K\eta^* m^*$ is well-defined.  For each $y \in Y_0$, let $\PP^y$ denote the Gaussian measure on $X$ with mean $m(y)$ and covariance operator $\hat K$.  In Theorem \ref{thm_mainresult}, we show that $\PP^y$ is a continuous disintegration of $\PP$ given $Y_0$.
	
	Suppose now that $M = \oo$, so that the operator $\eta^{-1}$ does not admit a continuous extension to all of $Y_0$.  Theorem 3.11 of \cite{tarieladze2007disintegration} gives the existence of a Borel-measurable linear operator $m$ and an operator $\hat K$ so that the Gaussian measure $\PP^y$ with mean $m(y)$ and $\hat K$ is a disintegration of $\PP$.  In Theorem \ref{thm_Minfinity}, we show that if there exists a continuous disintegration, then it must agree with $\PP^y$ on a set of full measure, and that the assumption $M=\oo$ implies that the conditional mean operator $m$ is discontinuous.  This results in a contradiction, thus $M < \oo$ is both a necessary and sufficient condition for there to exist a continuous disintegration of a Gaussian measure.
	

	\section{Probability Measures on Banach Spaces} \label{sect_probmeasinbanach}
	
	We now explore some of the general theory of Radon probability measures on Banach spaces.  Let $X$ be a Banach space, and let $\B(X)$ denote the Borel $\sigma$-algebra of $X$.  Continuous linear functionals of $X$ are measurable functions, hence random variables.  Let $\PP$ be a Radon probability measure on $X$ with the property that for all $f \in X^*$,
		\begin{equation} \label{weakordertwo}
			\EE|f|^2 = \int_X |f(x)|^2 \sD \PP(x) < \oo. \end{equation}
	This implies that every continuous linear functional has a finite variance and mean.  We recall that the \emph{support of the measure $\PP$} is the largest closed set in $X$ of full measure, and denote it by $\supp \PP$.
	
	\begin{env_thm}[Structure Theorem for Radon Probability Measures] \label{thm_structure}
		If $\PP$ is a Radon probability measure on $X$ which satisfies \eqref{weakordertwo}, then there exist an element $\mu \in X$ and a continuous linear operator $K : X^* \to X$ such that
			\begin{equation} \label{mean_K}
				f(\mu) = \EE(f) \qquad \mathrm{and} \qquad Kf = \int_X f(x) x \sD \PP(x) - f(\mu) \mu \end{equation}
		for all $f \in X^*$.  We call $\mu$ the \emph{mean} of $\PP$, and $K$ the \emph{covariance operator} of $\PP$.  It follows that 
			\begin{equation} \label{cov}
				f(Kg) = \EE(fg) - f(\mu)g(\mu) \end{equation}
		for all $f, g \in X^*$.    
		
		The space $\mu + KX^*$ is separable, and is dense in the support of $\PP$:
			\begin{equation} \label{suppeqn}
				\supp \PP \subseteq \mu + \overline{KX^*}. \end{equation}
		Consequently, $\PP(\mu + \overline{KX^*}) = 1$.  If $\PP$ is a Gaussian measure, then $\supp \PP = \mu + \overline{KX^*}$.
	\end{env_thm}
	\begin{proof}
		A measure $\PP$ which satisfies \eqref{weakordertwo} is called \emph{weak-order two}.  The existence of the mean vector $\mu$ is given by the Corollary in Section II.3.1 of \cite{vakhaniya1987probability}, and the existence of the covariance operator $K$ is given by Theorem 2.1 of Section III.2.1 of \cite{vakhaniya1987probability}.  The separability of the space $KX^*$ is Corollary 1 to that theorem.
		
		The statement about the support of a Gaussian measure is Theorem 1 of \cite{vakhania1975topological}.  The proof of \eqref{suppeqn} is part (a) of the proof of Theorem 1 of \cite{vakhania1975topological}.  The proof is simple and elegant so we reproduce it.  
		
		Without loss of generality, suppose that $\PP$ has mean zero.  Let $(KX^*)^\perp \subseteq X^*$ denote the annihilator of $KX^*$, defined below in \eqref{annihilator}.  If $f \in (KX^*)^\perp$, then $\int f(x)^2 \sD\PP(x) = f(Kf) = 0$, so $f(x) = 0$ for $\PP$-almost every $x \in X$.  Since the set $f^{-1}(0)$ is closed and has full measure, the support $\supp \PP$ is a subset of $f^{-1}(0)$.  Thus $f \in (\supp \PP)^\perp$.  \eqref{suppeqn} immediately follows.
	\end{proof}
	
	In addition to being a powerful technical result, the Structure Theorem presents a useful philosophy when working with Radon probability measures on Banach spaces:  many statements about probability can be reformulated in terms of the geometry of the linear space $\mu + \overline{KX^*}$.  This allows us to use linear algebra, functional analysis and, as we will see shortly, the theory of Hilbert spaces.
	
	When $\PP$ is a measure on the space of continuous functions with covariance function $c$ (e.g., Wiener measure, whence $c(t,s) = \min\{s,t\}$), the operator $K$ is the integral operator with kernel $c$.  This important special case is developed in Theorem \ref{thm_RF}.\newline

	For the remainder of this section, we assume that $\PP$ is a Radon probability measure on $X$ with mean zero and covariance operator $K$.

	Let $Y$ be a Banach space, and let $\B(Y)$ denote the Borel $\sigma$-algebra of $Y$.  Let $\eta : X \to Y$ be a continuous linear map from $X$ to $Y$.  Let $\PP_Y$ be the push-forward measure on $Y$ of $\PP$, defined by the equation
		$$\PP_Y(B) := \PP(\eta^{-1}(B))$$
	for every Borel set $B \in \B(Y)$.  This equation implies that the measure $\PP_Y$ satisfies the change of variable formula
		\begin{equation} \label{chvar}
			\int_{\eta^{-1}(B)} g(\eta x) \sD \PP(x) = \int_B g(y) \sD \PP_Y(y), \end{equation}
	for any integrable function $g : Y \to \R$.  Consequently, $\PP_Y$ has mean zero and covariance operator $\eta K \eta^*$.
	
	For a set $B$ of $X^*$, let 
		\begin{equation} \label{annihilator}
			B^\perp = \{ f \in X^* : f(Kg) = 0 \mathrm{~for~all~} g \in B \} \end{equation}
	be the \emph{annihilator} of $B$: the linear space of functionals uncorrelated with $B$,
		
	
	\begin{env_lem} \label{inj}
		When restricted to the subspace $KX^*$ of $X$, the map $\eta$ has kernel $K(\eta^* Y^*)^\perp$.  Consequently, on $K \eta^* Y^*$, $\eta$ is injective.  Define
			\begin{equation} \label{M} 
				M := \sup_{e \in Y^*} \left\{ \frac{\|K\eta^* e\|_X}{\|\eta K\eta^* e\|_Y} : e (\eta K\eta^* e) \ne 0 \right\}. \end{equation} 
		The inverse map $\eta^{-1} : \eta K\eta^* Y^* \to X$ has operator norm $M$.
	\end{env_lem}
	\begin{proof}
		Let $f \in X^*$.  For all $e \in Y^*$, 
			$$e(\eta Kf) = f(K\eta^* e)$$
		by the symmetry of the operator $K$, thus $f \in (\eta^* Y^*)^\perp$ exactly if $\eta(Kf) = 0$ in $Y$.  This proves that
			\begin{equation} \label{ker_eta}
				\ker \eta \cap KX^* = K(\eta^* Y^*)^\perp. \end{equation}
		
		The operator norm of the inverse map $\eta^{-1}$ on $\eta K \eta^* Y^*$ is given by
			$$\|\eta^{-1}\|_{\operatorname{op}} = \sup_{e \in Y^*} \left\{ \frac{\|K\eta^* e\|_X}{\|\eta K\eta^* e\|_Y} : \eta K\eta^* e \ne 0 \right\}.$$
		Let $M$ be as in \eqref{M}.  To see that $\|\eta^{-1}\|_{\operatorname{op}}$ and $M$ are equal, we apply the Schwarz inequality \cite{folland1999real} to the inner product on $Y^*$ generated by $\eta K \eta^*$:
			$$|e' \eta K \eta^* e|^2 \le |e' \eta K \eta^* e' | \, | e \eta K \eta^* e|.$$
		Thus, $\eta K \eta^* e \ne 0$ exactly if $e(\eta K \eta^* e) \ne 0$.
	\end{proof}
	
	
	Let 
		\begin{equation} \label{Y0}
			Y_0 = \overline{\eta K \eta^* Y^*}. \end{equation}
	Since the measure $\PP_Y$ has mean zero and covariance operator $\eta K \eta^*$, Theorem \ref{thm_structure} implies that it is supported on $Y_0$, so $\PP_Y(Y_0) = 1$.  Let $M$ be as in \eqref{M}, and suppose henceforth that 
		\begin{equation} \label{Mfinite}
			M < \oo. \end{equation}
	Define the linear map
		\begin{equation} \label{monY0}
			m : Y_0 \to X \end{equation}
	first by $m = \eta^{-1}$ on the dense subspace $\eta K \eta^* Y^*$ of $Y_0$, then extend continuously.  By Lemma \ref{inj}, the map $m$ is continuous with operator norm $M$.  Clearly, $\eta \circ m$ is the identity map on $Y_0$.  However, the map $m \circ \eta$ on $\overline{KX^*} \subseteq X$ is non-trivial.\newline



	The covariance operator $K$ defines a symmetric inner product $\langle f, g \rangle := f(Kg)$ on $X^*$.  This inner product is nonnegative-definite, and will be degenerate if the real-valued distribution of some $f \in X^*$ is atomic.  Nonetheless, it follows easily from the Schwarz inequality
		\begin{equation} \label{schwarz_inequality}
			|g(Kf)|^2 = |\langle g, f \rangle|^2 \le \langle g, g \rangle \, \langle f, f \rangle \end{equation}
	that $\langle f, f \rangle = 0$ if and only if $f \in \ker K$.  Thus the inner product is positive-definite on the quotient space $X^* / \ker K$.
	
	Let $H$ be the Hilbert-space completion of the inner product space $X^* / \ker K$, and let $\iota^* : X^* \to H$ be the inclusion map.  Define the unitary map $\iota : H \to X$ first on the dense subspace $\iota^* X^*$ by $\iota(\iota^* f) = Kf$, then extend it continuously to all of $H$.  The operator $K$ factors as $\iota \iota^*$.  We summarize this with the following commutative diagram:
	
		\begin{equation} \label{diagram}
			\begin{matrix} \xymatrix{Y^* \ar@{->}[r]^{\eta^*} & X^* \ar@{->}[rr]^{K} \ar@{->}[rd]^{\iota^*} && X\ar@{->}[r]^{\eta} & \ar@{->}@/^/[l]^m Y \\ && H \ar@{->}[ur]^{\iota} && } \end{matrix} \end{equation}
	where the domain of the map $m$ is $Y_0 \subseteq Y$.  
	
	The subspace $\iota H$ of $X$ is called the \emph{Cameron-Martin space} of $\PP$, and is a reproducing kernel Hilbert space \cite{berlinet2004rkh, janson1997ghs}.  The triplet $(\iota, H, X)$ is an \emph{abstract Wiener space} \cite{bell1987mc,gross2aws}.  Since $\iota H$ is dense in the separable Banach space $\overline{KX^*}$, the Hilbert space $H$ is separable.

	\begin{env_lem} \label{khat}
		The operator $\hat K : X^* \to X$ given by the formula
			\begin{equation} \label{Khatdef}
				\hat K = K - K \eta^* m^* \nonumber \end{equation}
		is well-defined.  Furthermore,
			\begin{equation} \label{meta}
				\hat K \le K \qquad \mathrm{and} \qquad m\eta K \eta^* m^* = K\eta^*m^*.\end{equation}
	\end{env_lem}
	The first statement of \eqref{meta} means that $f(\hat Kf) \le f(Kf)$ for all $f \in X^*$.
	
	\begin{proof}		
		Let $H_Y$ be the completion of $\iota^* \eta^* Y^*$ in $H$, and let $H_Y^\perp$ be its orthogonal complement.  Let $\pi : H \to H$ be the orthogonal projection map onto the subspace $H_Y$.  We claim that the two continuous maps $m \eta \iota$ and $\iota \pi$ from $H$ to $X$ are equal.  It suffices to check that that they are equal on the dense subspaces $\iota^* \eta^* Y^* \subseteq H_Y$ and $\iota^* (\eta^* Y^*)^\perp \subseteq H_Y^\perp$.  We calculate
			$$(m \eta \iota - \iota \pi)\iota^* \eta^* Y^* = m\eta K \eta^* Y^* - K \eta^* Y^* = 0$$
		since $\pi$ is the identity on $\iota^* \eta^* Y^*$ and $m \circ \eta$ is the identity on $K \eta^* Y^*$; and
			$$(m \eta \iota - \iota \pi)\iota^* (\eta^* Y^*)^\perp = m \eta K (\eta^* Y^*)^\perp - 0 = 0$$
		since $\pi$ kills $\iota^* (\eta^* Y^*)^\perp$ and $K (\eta^* Y^*)^\perp = \ker \eta \cap K X^*$ by Lemma \ref{inj}.  Thus 
			\begin{equation} \label{metaiota}
				m\eta\iota = \iota\pi \end{equation}
		on $H$.  By duality, the adjoint maps $\iota^* \eta^* m^*$ and $\pi \iota^*$ from $X^*$ to $H$ are also equal, so
			$$\hat K = K - K \eta^* m^* = K - \iota \iota^* \eta^* m^* = K - \iota \pi \iota^*$$
		is well-defined.  If we write $\pi^\perp : H \to H$ for the orthogonal projection map onto $H_Y^\perp$, then this shows that
			\begin{equation} \label{khatpiperp}
				\hat K = \iota \pi^\perp \iota^* \end{equation}
		since $K = \iota \iota^*$.  This representation implies that $\hat K \le K$.
		
		Finally, since $\pi^2 = \pi$,
			$$m\eta K \eta^* m^* = m\eta \iota \circ \iota^* \eta^* m^* = \iota \pi^2 \iota^* = \iota \pi \iota^* = K \eta^* m^*.$$
	\end{proof}
	
	Equation \eqref{khatpiperp} and Lemma \ref{inj} imply that
		\begin{equation} \label{kernel}
			\overline{\hat KX^*} = \ker \eta \cap \overline{K X^*}. \end{equation}
	Consequently,
		\begin{equation} \label{etamkhat}
			\eta \left( m(y) + \overline{\hat KX^*} \right) = y \end{equation}
	for all $y \in Y_0$, since $\eta \circ m$ is the identity on $Y_0$.\newline
	
	A Radon measure $\PP$ on $X$ is \emph{Gaussian} if every continuous linear functional $f \in X^*$ is a real-valued Gaussian random variable with respect to $\PP$.  Gaussian measures are completely described by their mean and covariance operators.  
	
	Not every continuous operator $\hat K : X^* \to X$ serves as the covariance operator for a Gaussian measure.  Nonetheless, the condition $\hat K \le K$ is sufficient, by Proposition 3.9 of \cite{tarieladze2007disintegration}.  We can now state and prove the main theorem of the paper.

	\begin{env_thmmain}[Existence of Continuous Disintegrations] \label{thm_mainresult}
		Let $X$ and $Y$ be Banach spaces, and let $\eta : X \to Y$ be a continuous linear map.  Let $\PP$ be a Gaussian probability measure on $X$ with mean zero and covariance operator $K$.  Suppose that
			\begin{equation} \label{M_thmmain}
				M = \sup_{e \in Y^*} \left\{ \frac{\|K\eta^* e\|_X}{\|\eta K\eta^* e\|_Y} : e (\eta K\eta^* e) \ne 0 \right\} < \oo. \end{equation}
		The operator $\eta^{-1} : \eta K \eta^* \to X$ is continuous with norm $M$.  Let $Y_0 := \overline{\eta K \eta^* Y^*}$, and let $m : Y_0 \to X$ denote the continuous extension of $\eta^{-1}$ to all of $Y_0$.  The map $\eta \circ m$ is the identity operator on $Y_0$.
		
		For any $y \in Y_0$, let $\PP^y$ be the Gaussian measure on $X$ with mean $m(y)$ and covariance operator $\hat K = K - K\eta^* m^*$.  The family of measures $\PP^y$ is a continuous disintegration given $Y_0$.  
	\end{env_thmmain}
	\begin{proof}
		Let $H$ be the Hilbert space described by the diagram \eqref{diagram} and let $\CC_H$ denote the cylinder algebra on $H$.  Let $\gamma : \CC_H \to [0,1]$ be the canonical Gaussian cylindrical measure on the Hilbert space $H$, i.e., with mean zero and covariance operator the identity $I$.  		
		
		Since the Cameron-Martin space $\iota H$ is dense in $\overline{KX^*}$, the push-forward cylindrical measure $\gamma \circ \iota^{-1}$ completely determines the measure $\PP$:
			\begin{equation} \label{PPgamma}
				\int_X f(x) \sD \PP(x) = \int_H f(\iota h) \sD \gamma(h) \end{equation}
		for any $f : X \to \R$ measurable with respect to the cylinder algebra $\CC$ of $X$.  If $f$ is a continuous linear functional, then the right side of \eqref{PPgamma} is further equal to $\int \langle \iota^* f, h \rangle \sD\gamma(h)$.
		
		Let $H_Y = \overline{\iota^* \eta^* Y^*}$ denote the subspace of $H$ generated by $Y$, and let $H_Y^\perp$ denote its orthogonal complement in $H$.  Let $\pi$ and $\pi^\perp$ be the orthogonal projection maps onto $H_Y$ and $H_Y^\perp$, respectively.  Define the push-forward cylindrical measures
			\begin{equation}
				\gamma_Y = \gamma \circ \pi^{-1} \qquad \mathrm{and} \qquad \gamma^\perp = \gamma \circ (\pi^\perp)^{-1} \end{equation}
		on $H$.  The cylindrical measure $\gamma_Y$ has mean zero, covariance operator $\pi$ and is supported on $H_Y$.  Similarly, $\gamma^\perp$ has mean zero, covariance operator $\pi^\perp$ and is supported on $H_Y^\perp$.
		
		We now exploit a fundamental fact of Gaussians:  orthogonality implies independence.  Let $k \in H_Y$ and $k' \in H_Y^\perp$.  The jointly Gaussian random variables $\langle k, \cdot \rangle$ and $\langle k', \cdot \rangle$ each have mean zero, and their covariance is $\langle k, I k' \rangle = 0$.  Consequently, the random variables $\langle k, \cdot \rangle$ and $\langle k', \cdot\rangle$ are independent.  Extending this analysis shows that for any $\CC_H$-measurable function $g : H \to \R$,
			\begin{equation} \label{gammaconv}
				\int_H g(h) \sD \gamma(h) = \int_H \int_H g(k+h) \sD \gamma^\perp(h) \D\gamma_Y(k). \end{equation}

		Clearly, $\PP_Y$ is the radonification of the push-forward cylindrical measure $\gamma_Y \circ (\eta \iota)^{-1}$ on $Y$.  Let $\PP^0$ be the radonification of $\gamma^\perp \circ \iota^{-1}$ on $X$.  This is the mean-zero Gaussian measure on $X$ with covariance operator $\hat K = \iota \pi^\perp \iota^*$, using the representation \eqref{khatpiperp}.  Define
			$$\PP^y(B) = \PP^0( m(y) + B )$$
		for any $B \in \B(X)$.  The measure $\PP^y$ is the Gaussian measure on $X$ with mean $m(y)$ and covariance operator $\hat K$.
		
		We now verify that $\PP^y$ is a disintegration.  By the structure theorem (Theorem \ref{thm_structure}), $\supp \PP^y = m(y) + \overline{\hat KX^*}.$  Thus by \eqref{etamkhat},
			$$\PP^y(\eta^{-1}(y)) \ge \PP^y\left( m(y) + \overline{\hat KX^*} \right) = 1,$$
		so $\PP^y$ is supported on the fiber $\eta^{-1}(y)$.
		
		The heart of the disintegration equation \eqref{disint} is the fact that $\gamma = \gamma^\perp * \gamma_Y$.  Let $f : X \to \R$ be measurable with respect to the cylinder algebra $\CC$ of $X$.  Then by equations \eqref{PPgamma} and \eqref{gammaconv}, 
			\begin{equation} \label{disint_calc1}
				\int_X f(x) \sD \PP(x) = \int_H f(\iota h) \sD \gamma(h) = \int_H \int_H f(\iota (k+h) ) \sD \gamma^\perp(h) \D \gamma_Y(k). \end{equation}
		For $\gamma_Y$-almost every $k$, $k = \pi k$.  We apply this to \eqref{disint_calc1}, as well as the identity \eqref{metaiota} that $\iota \pi = m \eta \iota$, to get
			\begin{equation} \label{disint_calc2}
				\int_H \int_H f(\iota \pi k + \iota h) \sD \gamma^\perp(h) \D\gamma_Y(k) = \int_H \int_H f(m\eta\iota k + \iota h) \sD \gamma^\perp(h) \D\gamma_Y(k). \end{equation}
		We now push forward to the Radon measures $\PP^0$ and $\PP_Y$, and use the definition $\PP^y = \PP^0( m(y) + \cdot)$, so that \eqref{disint_calc2} equals
			\begin{equation} \label{disint_calc3}
				\int_Y \int_X f(m(y) + x) \sD \PP^0(x) \D\PP_Y(y) = \int_Y \int_X f(x) \sD \PP^y(x) \D\PP_Y(y). \end{equation}
		Since the cylinder algebra $\CC$ generates the Borel $\sigma$-algebra $\B(X)$, this proves the disintegration equation \eqref{disint} for arbitrary integrable $f$.

		Finally, we show that $\PP^y$ satisfies the continuous disintegration property \eqref{contdis}.  Suppose $y_n \to y$ in $Y_0$.  The operator $m$ is continuous, so $m(y_n) \to m(y)$.  Let $f : X \to \R$ be a bounded, continuous function, so
			\begin{equation} \label{showingcontdisint}
				\lim_{n\to \oo} \int_X f(x) \sD \PP^{y_n}(x) = \lim_{n\to\oo} \int_X f(m(y_n) + x) \sD \PP^0(x) = \int_X f(m(y) + x) \sD \PP^0(x) = \int_X f(x) \sD \PP^y(x) \end{equation}
		by the bounded convergence theorem.  This proves that the measures $\PP^{y_n}$ converge weakly to $\PP^y$, which completes the proof.
	\end{proof}

	This theorem raises the natural question:  is $M < \oo$ a necessary condition for the existence of a continuous disintegration?  The next theorem demonstrates that for Gaussian measures, the answer is yes.  
	
	In the proof of Theorem \ref{thm_mainresult}, we used the fact that $M < \oo$ in order to define the conditional mean $m(y)$ and conditional covariance operator $\hat K$, as in \eqref{monY0} and \eqref{Khatdef}, respectively.  If we assume that $M = \oo$, then we must define the conditional mean and covariance using a different method.  The recent work \cite{tarieladze2007disintegration} of Tarieladze and Vakhania does exactly this, by working on the Hilbert space $H$ then using the identities $m\eta \iota = \iota \pi$ and $\hat K = \iota \pi^\perp \iota^*$ as the definitions of $m$ and $\hat K$.
	
	It is likely that the methods of Tarieladze and Vakhania can be adapted to a more general setting.  In that case, the Gaussian assumption can be weakened in the following theorem, though our proof does use the fact that the support of a mean-zero Gaussian measure $\PP$ is the entire linear space $\overline{KX^*}$, and not a proper subset.

	\begin{env_thmmain} \label{thm_Minfinity}
		Let $X$ and $Y$ be separable Banach spaces, and let $\eta : X \to Y$ a continuous linear map.  Let $\PP$ be a Gaussian measure on $X$ with mean zero and covariance operator $K$.  Suppose that
			\begin{equation}
				M = \sup_{e \in Y^*} \left\{ \frac{\|K\eta^* e\|_X}{\|\eta K\eta^* e\|_Y} : e (\eta K\eta^* e) \ne 0 \right\} = \oo. \end{equation}
		For any closed set $Y_0$ of full $\PP_Y$-measure, there does not exist a continuous disintegration $\PP^y$ on $Y_0$.
	\end{env_thmmain}
	\begin{proof}
		Let $Y_0$ be a closed subset of $Y$ of full $\PP_Y$-measure, and suppose $\tilde \PP^y$ is a continuous disintegration of $\PP$ on $Y_0$.

		The main result of \cite{tarieladze2007disintegration} is Theorem 3.11, which states that there exists a map $m : Y \to X$ and a covariance operator $\hat K$ such that the Gaussian measure $\PP^y$ with mean $m(y)$ and covariance $\hat K$ is a disintegration of $\PP$.  Furthermore, there exists a linear subspace $Y_1$ of $Y$ of full $\PP_Y$-measure such that the restriction of $m$ to $Y_1$ is Borel-measurable and $\eta(m(y)) = y$ for all $y \in Y_1$.
		
		Disintegrations are unique up to sets of measure zero \cite[Theorem 2.4]{tarieladze2007disintegration}, so there exists a set $Y_2$ of $Y$ of full $\PP_Y$-measure such that $\tilde \PP^y = \PP^y$ for all $y \in Y_2$.   Define the closed set
			$$Y' = \overline{Y_0 \cap Y_1 \cap Y_2},$$
		so that $\PP^y$ is a continuous disintegration of $\PP$ on $Y'$.  Since $Y'$ is a closed set of full $\PP_Y$-measure, it contains the linear space $\supp \PP_Y = \overline{\eta K \eta^* Y^*}$ as a subset.
		

		\begin{env_lem}
			There exists a sequence $y_n \in \eta K \eta^* Y^*$ such that $y_n \to 0$ but $\|m(y_n)\|_X \ge 1$ for all $n$.  Consequently, $m$ is discontinuous on $Y'$.
			
			Furthermore, the distance in $X$ from $m(y_n)$ to $\overline{\hat KX^*}$ is at least $1$ for all $n$.
		\end{env_lem}
		\begin{proof}
			As in the proof of Lemma \ref{khat}, let $H$ be the Hilbert space completion of the space $X^* / \ker K$ under the inner product generated by $K$, and let $\iota^* : X^* \to H$ be the inclusion map.  Define the unitary map $\iota : H \to X$ first on the dense subspace $\iota^* X^*$ by $\iota(\iota^* f) = Kf$, then extend it continuously to all of $H$.
		
			Let $H_Y$ be the completion of $\iota^* \eta^* Y^*$ in $H$.  Choose $e_i \in Y^*$ so that $h_i = \iota^* \eta^* e_i$ is an orthonormal basis in $H_Y$.  For all $y \in Y'$,
				$$m(y) = \sum_{i=1}^\oo e_i(y) \, K\eta^* e_i = \iota \left( \sum_{i=1}^\oo e_i(y) \, h_i \right);$$
			this follows from the proof of \cite[Theorem 3.11, Case 3]{tarieladze2007disintegration}.  If $\pi : H \to H$ is the orthogonal projection onto $H_Y$ in $H$, this formula implies that $m\eta \iota = \iota \pi$ on $H$.  Thus for $g \in Y^*$,
				\begin{equation} \label{munbounded}
					\|m(\eta K \eta^* g)\|_X = \| (m\eta\iota) \iota^* \eta^* g\|_X = \|\iota \pi \iota^* \eta^* g\|_X = \|K \eta^* g\|_X \end{equation}
			since $\pi$ is the identity on $\iota^* \eta^* Y^*$.
			
			Since $M = \oo$, there exist $g_n \in Y^*$ such that
				$$\|K\eta^* g_n \|_X \ge n \, \|\eta K \eta^* g_n \|_Y.$$
			Setting
				$$y_n = \frac 1 n \, \frac{\eta K \eta^* g_n}{\|\eta K \eta^* g_n\|}$$
			and applying \eqref{munbounded} completes the proof that $\|m(y_n)\| \ge 1$.
			
			It also follows from the proof of \cite[Theorem 3.11, Case 3]{tarieladze2007disintegration} that
				$$\hat K = K - \iota \pi \iota^*.$$
			Let $I$ denote the identity operator in $H$.  If $g \in Y^*$ and $f \in X^*$, then for $y = \eta K \eta^* g$,
				\begin{eqnarray*}
					\|m(y) - \hat K f\|_X^2 &=& \left\| \iota \iota^* \eta^* g - \left( \iota \iota^* f - \iota \pi \iota^* f \right) \right\|_X^2 \\ 
					&=& \| \iota^*\eta^*g - (I - \pi)\iota^* f\|_H^2 \\
					&=& \| \iota^*\eta^*g \|_H^2 + \|(I - \pi)\iota^* f\|_H^2 \\
					&\ge& \|m(y)\|_X^2 + 0,
				\end{eqnarray*}
			by the Pythagorean Theorem \cite{folland1999real}, since $I - \pi$ is the orthogonal projection onto $H_Y^\perp$.  Plugging in $y_n$ as above completes the proof of the second claim.
		\end{proof}
		
		Since $y_n \to 0$ in $\eta K \eta^* Y^* \subseteq Y'$ and $\PP^y$ is a continuous disintegration given $Y'$, $\PP^{y_n} \to \PP^0$ weakly.  By Theorem \ref{thm_structure}, $\PP^0$ is supported on $\overline{\hat K X^*}$.  Thus the open $\tfrac 1 2$-neighborhood of $\overline{\hat K X^*}$,
			$$U = \left\{ x \in X : \|x - x'\| < \tfrac{1}{2} \mathrm{~for~some~} x' \in \overline{\hat K X^*} \right\},$$ 
		has full $\PP^0$-measure.  Since $\PP^{y_n} \to \PP^0$ weakly and $U$ is open, $\liminf \PP^{y_n}(U) \ge \PP^0(U) = 1$.  However, $\PP^{y_n}$ is supported on $m(y_n) + \overline{\hat K X^*}$, which is distance at least $1$ from $\overline{\hat K X^*}$ by the preceding lemma, a contradiction.
	\end{proof}

	We now apply Theorem \ref{thm_mainresult} in the important context of random fields, which are stochastic processes defined on arbitrary parameter sets.  Let $T$ be a compact set.  Then a (continuous) random field $\xi_t$ on the parameter set $T$ is simply a random element of the Banach space $X = C(T, \R)$ of real-valued continuous functions on $T$.
	
	\begin{env_thmmain} \label{thm_RF} 
		Let $\xi_t$ be a Gaussian random field on a compact parameter set $T$ with mean zero and covariance function $c$:
			$$\EE \xi_t = 0 \qquad \mathrm{and} \qquad \EE \xi_t \xi_s = c(t,s).$$
		Suppose that $\xi_t$ is almost-surely continuous.  Let $S$ be a closed subset of $T$, and suppose that
			\begin{equation} \label{RF_M}
				M = \sup_{s \in S} \frac{\sup_{t \in T} |c(s,t)|}{\sup_{s' \in S} |c(s,s')|} < \oo. \end{equation}
		There exists a closed family $Y_0$ of functions such that with probability one, $\xi|_{S} \in Y_0$, and the regular conditional probability $\PP^y := \PP\big( \cdot \big|\, \xi|_{S} = y \big)$ is a well-defined Gaussian measure which varies (weakly) continuously in $y \in Y_0$. 
				
		Furthermore, there exists a function $m(y,t)$, linear in $y$ and jointly continuous in $y$ and $t$, such that
			\begin{equation} \label{RF_mean}
				 m(y,t) = \EE\big(\xi_t \big|\, \xi|_S = y\big) \end{equation}
		for all $t \in T$, and there exists a covariance function $\hat c$ (independent of $y$) such that
			\begin{equation} \label{RF_cov}
				\hat c(t,s) = \EE\big(\xi_t \xi_s \big|\, \xi|_{S} = y\big) - m(y,t)m(y,s) \end{equation}
		and $\hat c(t,s) \le c(t,s)$ for all $t,s \in T$.  The function $m(y,\cdot)$ is a bounded extension of $y$, in the sense that
			\begin{equation} \label{RF_boundedmean}
				m(y,s) = y(s) \mbox{~for all $s \in S$,} \qquad \mathrm{and} \qquad \sup_{t \in T} |m(y,t)| \le M \sup_{s \in S} |y(s)|. \end{equation}
	\end{env_thmmain}
	\begin{proof}
		Consider the Banach space $X = C(T, \R)$.  By the Riesz representation theorem \cite{folland1999real}, the dual space $X^*$ has a representation as the space of Radon measures on $T$, so for all $f \in X^*$ there exists a Radon measure $\lambda_f$ so that $f(x) = \int_A x(t) \sD \lambda_f(t)$.  Define the operator $K : X^* \to X$ by
			$$(Kf)(t) = \int_A c(t,s) \sD \lambda_f(s).$$
		That is, $K$ is the integral operator with kernel $c$.  Consequently, the measure $\PP$ on $X$ is a Radon probability measure with mean zero and covariance operator $K$.  Let $\delta_t$ represent the evaluation functional, defined by $\delta_t x = x(t)$; equivalently, $\delta_t$ represents the Dirac point-mass measure with an atom at $t \in T$.  Thus $c(t,s) = \delta_t(K \delta_s)$.  
		
		
		Let $Y = C(S, \R)$, and let $\eta : X \to Y$ be the restriction map, defined by $(\eta x)(s) = x(s)$ for all $x \in X$ and $s \in S$.  Let $\PP_Y = \PP \circ \eta^{-1}$ denote the push-forward of $\PP$ onto $Y$.  Using the same notation as above, denote the evaluation functionals on $Y$ by $\delta_s$.  Let $M$ be defined by \eqref{M_thmmain}.  Since the linear spans of $\{\delta_t\}_{t\in T}$ and $\{\delta_s\}_{s\in S}$ are dense in $X^*$ and $Y^*$ \cite{reed-methods}, respectively, $M$ takes the form \eqref{RF_M}.  
		
		In this context, the space $Y_0 = \overline{\eta K \eta^* Y^*}$ takes the form
			\begin{equation} \label{RF_Y0}
				Y_0 = \overline{\span\{c(s,\cdot)\}} \subseteq Y, \end{equation}
		where the span is over $s \in S$.  The space $Y_0$ has full $\PP_Y$-measure.
		
		By assumption, $M < \oo$, so Theorem \ref{thm_mainresult} applies.  Thus there exists a continuous disintegration $\PP^y$ on $Y_0$, and there exist continuous linear operators $m : Y_0 \to X$ and $\hat K : X^* \to X$ so that for all $y \in Y_0$, the measure $\PP^y$ has mean $m(y)$ and covariance operator $\hat K$.  Define the functions $m(y,t) := \delta_t (m(y))$ and $\hat c(t,s) := \delta_t(\hat K \delta_s)$.  Since $\eta \circ m$ is the identity operator on $Y_0$, and $m$ has operator norm $M$, the statements \eqref{RF_boundedmean} immediately follow.
	\end{proof}
	
	Suppose that $T$ is a subset of an abelian group.  We say that a random field $\xi_t$ is stationary if its covariance function satisfies
		$$c(t,s) = c(t + z, s+z),$$
	whenever $t$, $s$, $t+z$ and $s+z$ all belong to $T$.

	\begin{env_cormain} \label{cor_stationary}
		If $\xi_t$ is a stationary Gaussian random field, then $M = 1$, so the above theorem applies to $\xi_t$.
	\end{env_cormain}
	\begin{proof}
		Since the covariance operator $K$ defines an inner product, the Schwarz inequality \eqref{schwarz_inequality} implies
			\begin{equation}
				|c(s,t)|^2 \le c(s,s) \, c(t,t). \end{equation}
		If the field is stationary, then $c(t,t) = c(s,s)$.  Consequently, for each $s \in S$,
			$$\sup_{t \in T} |c(s,t)| = |c(s,s)| \qquad \mathrm{and} \qquad \sup_{s' \in S} |c(s,s')| = |c(s,s)|.$$
		The ratio of these two quantities is always equal to $1$, so $M = 1$.
	\end{proof}

		\textbf{Acknowledgements.}  
		
		The author thanks Janek Wehr and Joe Watkins for many useful discussions and helpful feedback, and Mark Meckes \cite{meckesMO} for the simple proof of equation \eqref{showingcontdisint}.  The author is particularly indebted to Nicholas Vakhania and Vaja Tarieladze for their fundamental work on the subject of Radon probability measures on Banach spaces, and he hopes that this article does justice to their ideas.
		
		The author was supported by NSF VIGRE Grant No. DMS-06-02173 at the University of Arizona, and by NSF PIRE Grant No. OISE-07-30136 at the Courant Institute of Mathematical Sciences.

\end{document}